\documentclass[12pt]{article}
\usepackage{amsfonts,amssymb}

\newtheorem{theo}{Theorem}[section]

\newtheorem{lem}[theo]{Lemma}
\newtheorem{cor}[theo]{Corollary}
\newtheorem{defn}[theo]{Definition}

\newcommand{\norm}[1]{\Vert#1\Vert}
\newcommand{\alg}{{\bf g}}
\newcommand{\ala}{{\bf a}}
\newcommand{\aln}{{\bf n}}
\newcommand{\cross}{\ltimes}
\newcommand{\diam}{{\rm diam}}
\newcommand{\lin}{{\rm lin}}
\newcommand{\ad}{{\rm ad}}
\newcommand{\Ad}{{\rm Ad}}
\title{Calder\'on Zygmund decompositions on amenable groups
}
\author{
W. Hebisch\thanks{Partially supported by KBN
grant {\tt 2 P03A 058 14} and European Commission via TMR network
``Harmonic analysis''}\\
{\small  Institute of Mathematics, Wroc\l aw University, Wroc\l aw, Poland}
}
\begin{document}
\maketitle

\begin{abstract} 
We propose a simple abstract version of Calder\'on--Zygmund theory,
which is applicable to spaces with exponential volume growth, and then
show that amenable Lie groups can be treated within this
framework.
\end{abstract}

\section{Abstract Calder\'on--Zygmund theory}
\begin{defn}\label{CZ}
We say that the space~$M$ with metric~$d$ and Borel measure~$\mu$ has
the \emph{Calder\'on--Zygmund property} if there exists a constant~$C$
such that for every~$f$ in~$L^1$ and for every
$\lambda>C\frac{\|f\|_{L^1}}{\mu(M)}$ ($\lambda>0$ if $\mu(M)=\infty$)
we have a decomposition $f=\sum f_i+g$, such that there exist sets
$Q_i$, numbers $r_i$, and points $x_i$ satisfying:
\begin{itemize}
\item $f_i=0$ outside $Q_i$,
\item $\int f_i\,d\mu=0$,
\item $Q_i\subset B(x_i,Cr_i)$,
\item $\sum \mu(Q^{*}_i) \leq C\frac{\|f\|_{L^1}}{\lambda}$, where
$Q^{*}_i=\{x: d(x,Q_i)< r_i\}$,
\item $\sum \|f_i\|_{L^1}\leq C\|f\|_{L^1}$,
\item $|g|\leq C\lambda$.
\end{itemize}
Since $g=f-\sum f_i$, we have $\norm{g}_{L^1}\leq C'\norm{f}_{L^1}$,
hence $\norm{g}_{L^2}^2\leq C''\lambda\norm{f}_{L^1}$.
\end{defn}

\begin{theo}\label{amencz}
If $G$ is a connected amenable Lie group,
$d$ is a right-invariant optimal control metric on $G$ then $G$ 
(with $d$ and left Haar measure) satisfies Calder\'on--Zygmund property.
\end{theo}

From now to the end of the section we consider a fixed separable metric
space $M$ with a metric $d$ and measure $\mu$. We assume that all balls
have finite measure. All operations on sets are meant as operations in
the lattice of measurable sets on $M$.

%\begin{defn}
We say that a set $Q$ is a {\it doubling set} with a constant $C$
iff $\mu(\{x: d(x, Q) \leq {1\over C}\diam(Q)\}) \leq C\mu(Q)$.

We say that a family ${\cal A}$ of sets is a {\it doubling family}
with a constant $C$ if for each $Q \in {\cal A}$ it satisfies the
following two conditions:
\begin{itemize}
\item the set $$\tilde Q = 
\bigcup_{R\in {\cal A}, R \cap Q \ne \emptyset, \mu(R) \leq 2\mu(Q)} R$$
is contained in a set $S\in{\cal A}$ such that $\mu(S) \leq C\mu(Q)$.
\item either $\mu(M)\leq C\mu(Q)$ and $M \in {\cal A}$ or
there is $R\in {\cal A}$ such that
$Q\subset R$ and $2\mu(Q) \leq \mu(R) \leq C\mu(Q)$.
\end{itemize}

{\bf Remark} We get equivalent definition of doubling family
if in the definition of ${\tilde Q}$ we replace $\mu(R) \leq 2\mu(Q)$
by $\mu(R) \leq \mu(Q)$.  Namely, choose $S\in {\cal A}$ 
using second point of the definition
so $Q \subset S$ and
% Enlarging $S$ if needed using second point of the definition we may
% assume that 
$\min(\mu(M), 2\mu(Q)) \leq \mu(S) \leq C\mu(Q)$.
Now, if $R\cap Q\ne \emptyset$,
then $R\cap S\ne \emptyset$.  Also, if $\mu(R) \leq 2\mu(Q)$, then
$\mu(R) \leq \mu(S)$.  So,
${\tilde Q}_{\rm def\ 1} \subset {\tilde S}_{\rm def\ 2}$, where
${\tilde Q}_{\rm def\ 1}$ uses first definition while ${\tilde S}_{\rm def\ 2}$
uses condition $\mu(R) \leq \mu(Q)$.  So, applying weaker version of first
condition to $S$ we get stronger version for $Q$ with $C$
replaced by $C^2$.

We say that a family ${\cal A}$ of sets is {\it dense} if for a.e $x$ and
all $r>0$ there is a $Q\in{\cal A}$ such that $x\in Q$ and $\diam(Q) \leq r$.

%\end{defn}

\begin{theo}
Assume ${\cal A}$ is a doubling family. The maximal function
$$
M_{\cal A}f(x) = \sup_{Q\in{\cal A}, x \in Q} {1\over \mu(Q)} \int_Q |f|$$
is of weak type $(1, 1)$.

If ${\cal A}$ is dense then for each $f\in L^1(\mu)$
$$
\lim_{r \rightarrow 0+} \sup_{Q\in{\cal A}, \diam{Q}\leq r, x \in Q}
{1\over \mu(Q)} \int_Q |f| = |f|(x)$$
except possibly for a set of $\mu$ measure $0$.

If ${\cal A}$ is dense and consists of doubling sets (with a common
doubling constant $C$) then $M$ satisfies Calder\'on--Zygmund property.
\end{theo}

{\it Proof}: Let ${\cal S}_0$ be the set of all $R\in {\cal A}$ such
that 
$$\int_R |f| > \lambda \mu(R).$$
Let $v_0 = \sup_{R\in {\cal S}_0} \mu(R)$. $v_0$ is finite since
for all $R\in {\cal S}_0$ $\mu(R) < \|f\|_{L^1}/\lambda$. We take as
$R_0$ an element of ${\cal S}_0$ such that $\mu(R_0)>v_0/2$ (we choose
arbitrarily among all $R$ satisfying this bound). Now, we proceed
inductively: ${\cal S}_{i}$ is the set of all $R\in {\cal S}_{i-1}$ such
that $R$ is disjoint with all $R_j$, $j=0,\dots,i - 1$. 
Put $v_i = \sup_{R\in {\cal S}_i} \mu(R)$ and choose $R_i$ such
that $\mu(R_i)>v_i/2$. 

We choose as $Q_i$ element of ${\cal A}$ such that $\tilde R_i
\subset Q_i$ and $\mu(Q_i) \leq C\mu(R_i)$ (possible since our family
is a doubling family).

$R_i$ are disjoint, so
$$
\sum_i \mu(R_i) \leq {1\over \lambda}\sum_i \int_{R_i}|f| 
\leq {1\over \lambda} \|f\|_{L^1}$$
and
$$
\sum_i \mu(\tilde R_i) \leq \sum_i \mu(Q_i) \leq C\sum_i \mu(R_i) 
\leq {C\over \lambda} \|f\|_{L^1}$$

If $\int_R |f| > \lambda \mu(R)$ then by construction $R$ intersects
some $R_i$ such that $\mu(R) \leq 2\mu(R_i)$, so 
$R \subset \tilde R_i \subset Q_i$.
Hence, putting $E = \bigcup_i \tilde R_i$ we have $M_{\cal A}f \leq \lambda$
outside $E$ and $\mu(E) \leq {C\over \lambda} \|f\|_{L^1}$ which
gives the first claim.

%First, note the following Wiener type lemma:

%\begin{lem}\label{winerl}
%If the set $E$ is a sum of a collection ${\cal S} \subset {\cal A}$,
%$\mu(E) > t$, then there exists
%a finite subcollection $Q_j\in {\cal S}$, $j=1,\dots,n$ such that
%$Q_j$ are disjoint, $\mu(E\cap \bigcup_{j=1}^n \tilde Q_j) > t$,
%and $\sum_{j=1}^{n}\mu(Q_j)>t/C$.
%\end{lem}

%The \ref{winerl} immediately implies that $M_{\cal A}$ is of weak
%type $(1, 1)$.

If the family ${\cal A}$ is dense then the second claim is valid for
continuous functions with bounded support. Since such functions are
dense in $L^1(\mu)$ we get the second claim using weak type $(1, 1)$ of
$M_{\cal A}$.

We put $U_0 = \tilde R_0$ and
$U_i = \tilde R_i - \bigcup_{j<i} \tilde R_j$. Let $h_i(x) = f(x)$ 
on $U_i$ and $0$ otherwise. We claim that
$$
\int_{Q_i} |h_i| \leq C\lambda\mu(Q_i)$$
We need to consider three cases. In the first case $C\mu(Q_i) > \mu(M)$ and
$$
\int_{Q_i} |f| > C\lambda\mu(Q_i).$$
Then
$$
\|f\|_{L^1} > C\lambda\mu(Q_i) > \lambda\mu(M)$$
and there is no need to do Calder\'on-Zygmund decomposition
($\lambda$ is out of range).
In the second case $\mu(Q_i) > 2\mu(R_i)$ and we 
put $Q = Q_i$.  In the third case $\mu(Q_i) \leq 2\mu(R_i)$
we take as $Q$
an element of ${\cal A}$ such that $Q_i\subset Q$, $\mu(Q) < C\mu(Q_i)$
and $\mu(Q) \geq 2\mu(Q_i) \geq 2\mu(R_i)$. If
$$
\int_{Q} |f| \leq \lambda\mu(Q),$$
then
$$
\int_{Q_i} |h_i| \leq \int_{Q} |f| 
\leq \lambda\mu(Q) \leq C\lambda\mu(Q_i).$$
If the inequality above does not hold, then
$$
\int_{Q} |f| > \lambda\mu(Q)$$
and by our construction $Q$ intersects some of $R_j$ with $j<i$ and
$2\mu(R_j) > \mu(Q)$. So
$Q_i \subset Q \subset \tilde R_j$
and in this case $h_i = 0$ (since
$U_i = \emptyset$).

Next, let $\chi_{R_i}$ be the indicator (characteristic) function
of $R_i$. We put
$$
f_i = h_i - {\int h_i \over \mu(R_i)}\chi_{R_i},$$
$$
g = f - \sum f_i,$$
$$r_i = {1\over C}\diam(Q_i).$$
For $x_i$ we choose an arbitrary point from $Q_i$, which finishes
construction of Calder\'on-Zygmund decomposition. 

We need to check that conditions of Calder\'on-Zygmund decomposition
are satisfied. 
By construction $f_i = 0$ outside $Q_i$ and $\int f_i = 0$. Also
$Q_i \subset B(x_i, Cr_i)$. Since each $Q_i$ is doubling (with constant C)
we have
$$\mu(Q_i^{\star}) \leq C\mu(Q_i) \leq C^2\mu(R_i)$$
and
$$\sum_i \mu(Q_i^{\star}) \leq C^2\mu(R_i) 
\leq {C^2 \over \lambda}\|f\|_{L^1}.$$
Also
$$
\|f_i\|_{L^1} \leq 2\|h_i\|_{L^1} = 2\int_{U_i} |f|$$
so (noting that $U_i$ are disjoint)
$$
\sum_i \|f_i\|_{L^1} \leq 2\sum_i \int_{U_i} |f| \leq 2\|f\|_{L^1}.$$
For $x$ outside $E$ we have
$|g|(x) \leq |f|(x) \leq M_{\cal A}f(x) \leq \lambda$ (the second inequality
follows from second claim of our theorem). On $E$, since
$R_i$ are disjoint, we have
$$
\sup_{x \in E}|g|(x) \leq \sup_i |{\int h_i \over \mu(R_i)}|
\leq \sup_i {C\lambda\mu(Q_i) \over \mu(R_i)}$$
$$
\leq C^2\lambda.$$
$\diamond$

\section{Base case}

In this section we will prove a special case of the main theorem.

\begin{theo}\label{basecase}
Assume that  $G$ is a connected and simply connected solvable Lie group,
$N$ and $W$ are connected and
simply connected nilpotent subgroups of $G$, $N$ is a normal subgroup,
$G=WN$. Also assume that $G$ and $N$ are equipped with right
invariant Riemannian 
metrics $d_G$ and $d_N$ such that for $x, y\in N$
$$
\exp(C_1d_G(x,y)) \leq 1 + d_N(x,y)  \leq \exp(C_2(d_G(x,y) + 1))
$$
Then $G$ contains dense doubling family of doubling sets
and $G$ satisfies Calder\'on--Zygmund property.
\end{theo}

Before we prove \ref{basecase} we need some preparations.
Let $W_0 = W/(W\cap N) = G/N$ and let $\pi:G\mapsto W_0$ be the quotient
mapping. On $W_0$ we put quotient metric (from $W$):
$d_{W_0}(x,y) = \inf_{v, z \in W, \pi(v) = x, \pi(z) = y} d_G(v, z)$.

The following lemma is known, but we provide proof for
convenience.
% (Kuratowki, Wstep, Chapter XVI, exercises 34 and 41).
\begin{lem}\label{messel0}
Let $X$ and $Y$ be complete separable metric spaces, $Y$ compact,
$R \subset X \times Y$ a compact subset.  If for each $y\in Y$ the set
$R_y = \{x : (x, y) \in R \}$ is nonempty,
then there exist Borel measurable
$\chi: Y \rightarrow X$ such that $(\chi(y), y) \in R$ for all $y \in Y$.
\end{lem}
{\it Proof}:  Without loss of generality we may assume that
$X = H = [0,1]^\infty$.  Namely, it is well known that any
complete separable metric space is homeomorphic to a 
subset of $H$, so we can treat $X$ as a subset of $H$.
Then image of $R$ is a compact subset of $H\times Y$.
Now, again without loss of generality
we may assume that $X = [0,1]$.  Namely, it is well known that
there is continuous function $f$ from $[0,1]$ onto $H$.  Putting
$h(x, y) = (f(x), y)$ we get mapping from $[0,1]\times Y$ onto
$H\times Y$.  $S = h^{-1}(R)$ is a closed (hence compact)
subset of $[0,1]\times Y$.  Once we build
$\chi$ for $S$ the composition $f\circ\chi$ gives result for $R$.
So now $X = [0,1]$.  We claim that $\chi(y) = \inf(R_y)$ is
Borel measurable and has required properties.  First,
since $R_y$ is nonempty $\chi$ is well defined.  Next,
since $R_y$ is compact, we have $\inf(R_y) \in R_y$ so
indeed $(\chi(y), y) \in R$.
To show that $\chi$ is Borel measurable we need extra
construction.  For each rational $q \in [0, 1]$ consider
set $T_q = R \cap ([0, q]\times Y)$ and let $F_q$ be projection
of $T_q$ on $Y$.  We put $\phi_q(y) = q$ for $y \in T_q$ and
$\phi_q(y) = 1$ otherwise.  Since $T_q$ and $F_q$ are compact
$\phi_q$ is Borel measurable.  We claim that
$\chi(y) = \inf_{q}\phi_q(y)$.  Clearly $(\chi(y), y) \in R$,
so for any $q \geq \chi(y)$ we have $(\chi(y), y) \in T_q$,
so $y\in F_q$ and $\phi_q(y) = q$.
Since $\chi(y) = \inf(R_y)$, for $q < \chi(y)$ we have
$[0, q] \cap R_y = \emptyset$, so $[0, q] \times \{y\} \cap R = \emptyset$.
Hence $y \notin F_q$
and $\phi_q(y) = 1$.  Together,
$$
\inf_q \phi_q(y) = \inf_{q \geq \chi(y)} q = \chi(y).$$
Since infimum of countable
family of Borel measurable functions is Borel measurable
$\chi$ is Borel measurable.
$\diamond$

\begin{lem}\label{messel}
There exists a Borel measurable map $\chi$ from $W_0$ into $W$ such
that for all $x\in W_0$, $\pi\circ\chi(x)=x$ 
and $d_G(\chi(x), e) = d_{W_0}(x, e)$.
\end{lem}
{\it Proof}: Consider the relation 
$R(x, y) \iff (x = \pi(y) \wedge d_G(y, e) = d_{W_0}(x, e))$.
It is easy to see that $R$ is closed, and that for each $x$ the set 
$\{y : R(x, y)\}$ is a nonempty compact. Also, we may restrict
$x$ to stay in a compact set: we simply cover $W_0$ by a countable family
of compact sets $K_n$, build $\chi_n$ in each $K_n$ separately and
then glue them in Borel measurable way (for example by taking 
smallest $n$ such that $x$ is in domain of $\chi_n$).
Now the claim follows from \ref{messel0}.
$\diamond$

\begin{lem}\label{ballprod}
$B_G(e, r) \subset \chi(B_{W_0}(e, r))(B_G(e, 2r)\cap N)$.
\end{lem}

{\it Proof}: If $x \in B_G(e, r)$, then $\pi(x) \in B_{W_0}(e, r)$.
Let $z = \chi(\pi(x))$. We have $x = z(z^{-1}x)$.
Next $\pi(z^{-1}x) = \pi(z)^{-1}\pi(x) = \pi(x)^{-1}\pi(x) = e$, so
$z^{-1}x \in N$. By \ref{messel}, $d_G(z, e) = d_G(x, e)$, so
$d_G(z^{-1}x, e) \leq d_G(z, e) + d_G(x, e) = 2d_G(x, e)$.
$\diamond$

\begin{lem}
We can normalize left invariant Haar measures $\mu$ ($\lambda$, $\nu$)
on $G$ ($W_0$, $N$ respectively) in such a way that the following holds:
If $\chi: W_0 \mapsto W$ is a Borel measurable function such
that $\pi\circ\chi(x) = x$, than
$$\mu(\chi(Q)R) = \lambda(Q)\nu(R)$$
for all Borel measurable $Q \subset W_0$, $R\subset N$.
The same normalization works for all such $\chi$.
\end{lem}

{\it Proof}:  Fix $\lambda$ and $\nu$.  Note that mapping
$(w, n) \mapsto \chi(w)n$ is a Borel measurable one-to-one
mapping from $W_0\times N$ onto $G$ with Borel measurable
inverse.
So formula $\mu(\chi(Q)R) = \lambda(Q)\nu(R)$
defines measure $\mu$ as a transport of $\lambda\otimes\nu$.
We need to check that resulting $\mu$ is left invariant on $G$
and does not depend on choice of $\chi$.  By definition and
Fubini theorem, for nonnegative Borel measurable $f$ on $G$
we have
$$
\int_G fd\mu = \int_{W_0} \int_N f(\chi(w)n)d\nu(n)d\lambda(w)
$$
If $h : W_0 \rightarrow N$ is Borel measurable by
left invariance of $\nu$ we have
$$
\int_{W_0} \int_N f(\chi(w)n)d\nu(n)d\lambda(w)
= \int_{W_0} \int_N f(\chi(w)h(w)n)d\nu(n)d\lambda(w).
$$
If $\chi_1$ and $\chi_2$ are two different choices for $\chi$
we have $\chi_2(w) = \chi_1(w)h(w)$ with $h(w) \in N$
so by the formula above $\mu$ does not depend on choice of $\chi$.
Similarly, if $z \in N$ putting $h(w) = \chi(w)^{-1}z^{-1}\chi(w)$ we have
$\pi(h(w)) = \pi(\chi(w))^{-1}\pi(z^{-1})\pi(\chi(w)) = e$, so
$h(w) \in N$.  Then $z\chi(w)h(w)n = \chi(w)n$ and by the formula above
$$
\int_G f(zx)d\mu(x) =
\int_{W_0} \int_N f(z\chi(w)n)d\nu(n)d\lambda(w)
$$
$$
= \int_{W_0} \int_N f(\chi(w)n)d\nu(n)d\lambda(w) = \int_G fd\mu
$$
so measure $\mu$ is
left invariant under action of $N$.

Finally, if $y \in W$ then $y\chi(w) = \chi(\pi(y)w)h(w)$.  Again
$$
\pi(h(w)) = (\pi(\chi(\pi(y)w))^{-1}\pi(y\chi(w)) =
(\pi(y)w)^{-1}\pi(y)w = e
$$
so $h(w) \in N$ and by the
formula above and left invariance of $\lambda$ we have
$$
\int f(yx) d\mu(x) =
\int_{W_0} \int_N f(\chi(\pi(y)w)h(w)n)d\nu(n)d\lambda(w)
$$
$$
= \int_{W_0} \int_N f(\chi(w)n)d\nu(n)d\lambda(w) = \int f(x)d\mu(x)
$$
so $\mu$ is left invariant under action of $W$.  Since
$G=WN$ this means that $\mu$ is left invariant under action of
$G$.
$\diamond$

\begin{lem}\label{unidouble}
Consider a nilpotent Lie group equipped with an optimal control
metric $d$.
The doubling property:
$$
|B(x,2r)| \leq C|B(x,r)|$$
holds with a constant $C$ which depends only on the dimension of the group.
In particular, the constant is independent of $d$.
\end{lem}
{\it Proof}: This follows via transference from free nilpotent group.
Namely, let $G$ be a nilpotent Lie group of dimension $n$ and let
$F$ be free nilpotent group of step $n$ on $n$ generators.  Let $d$
be associated to vector fields $X_1,\dots,X_n$.  We map generators of
$F$ to $X_1,\dots,X_n$.  Then we get homomorphizm $\pi$ from $F$ to $G$,
such that $B_{G}(r, e_G) = \pi(B_F(r, e_F))$ (where we use subscripts
to distinguish objects in $F$ and $G$).  We may consider $G$ as an
$F$-space, where $x \in F$ acts on $G$ by left multiplication by $\pi(x)$.
By Lemma 1.1 in \cite{Gui} the following inequality holds, whenever
locally compact group $F$ acts on locally compact space $G$ equipped
with $F$ invariant measure and $A, B\subset F$, $Y\subset G$:
$$
|A||BY| \leq |BA||A^{-1}Y|.$$
Putting $A = A^{-1} = B_F(r, e_F)$, $B=B_F(2r, e_F)$, $Y= \{e_G\}$
we get:
$$
|B_F(r, e_F)||B_G(2r, e_G)| = |A||BY| \leq |BA||A^{-1}Y| = 
|B(3r, e_F)||B_G(r, e_G)|$$
so
$$
|B_G(2r, e_G)|  \leq C|B_G(r, e_G)|$$
where $C = {|B(3r, e_F)| \over |B_F(r, e_F)|}$.
$\diamond$

We say that a finite sequence of metrics $d_i$, $i=0,\dots, k$ on $N$
is a {\it doubling chain of metrics} iff
$$2d_{i+1}(x, y) \leq d_{i}(x, y) \leq 16d_{i+1}(x, y)$$
for all $x$ and $y$.

\begin{lem}\label{doublechain}
Let $d$ and $\rho$ be two right-invariant Riemannian metrics on $N$.
If there is $m \geq 2$ such that
$ m^2\rho \geq d \geq m \rho $,
then there exists
doubling chain of metrics such that $d_0 = d$ and $d_k = \rho$.
% and
% that Riemannian measure of $d_i$ is at least twice as big as 
% Riemannian measure of $d_{i+2}$.
\end{lem}
{\it Proof}: Right-invariant Riemannian metric on a group
is uniquely determined by
corresponding scalar product at $e$.  So it is enough to prove
analogous lemma for quadratic forms.  But two positive quadratic forms
can be diagonalized simultaneously (we diagonalizes one form first, and
then use orthogonal transformations with respect to this form to
diagonalize the other one).  For quadratic forms in diagonal form
construction of required chain is straightforward.
$\diamond$

\begin{lem}\label{ad-metric}
If $x\in W$, and $d_1(n_1, n_2) = d_N(x n_1 x^{-1}, x n_2 x^{-1})$,
then there is $C$ such that $d_1 \leq \exp(Cd_G(x, e))d_N$
\end{lem}
{\it Proof}:  This is well known.
%Since $d_1$ and $d_N$ are determined by the scalar
%product at $e$ it is enough to prove similar bound for 
$\diamond$

The following theorem is due to M. Christ (\cite{ChT} Theorem 11):

\begin{theo}
Let $X$ be a space of homogeneous type.  There exists a collection
of open subsets $\{Q_{\alpha,k} : k \in {\cal Z}, \alpha \in I_k\}$,
and constants $\delta \in (0, 1)$, $a_0 > 0$ and
$C< \infty$, such that
\begin{itemize}
\item $\forall_k \mu(X - \bigcup_{\alpha} Q_{\alpha,k}) = 0$.
\item If $l > k$, then either $Q_{\beta,l} \subset Q_{\alpha,k}$
 or $Q_{\beta,l} \cap Q_{\alpha,k} = \emptyset$.
\item For each $(k,\alpha)$ and each $l < k$, there is unique
 $\beta$ such that $Q_{\alpha,k} \subset Q_{\beta,l}$.
\item $\diam(Q_{\alpha,k})\leq C\delta^k$.
\item Each $Q_{\alpha,k}$ contains some ball $B(z_{\alpha,k}, a_0\delta^k)$.
\end{itemize}
\end{theo}

We say that $\{Q_{\alpha,k} \}$ is a family of dyadic cubes.  For
each $Q \in \{Q_{\alpha,k}\}$ we choose a point $x_Q \in Q$ which
we will call centerpoint of $Q$.

{\bf Remark} First and second condition together mean that 
$\mu(Q_{\alpha,k} - \bigcup_{\beta\in J_\alpha}Q_{\beta,k+1}) = 0$
where $J_\alpha = \{\beta: Q_{\beta,k+1} \subset Q_{\alpha,k}\}$.

In the sequel we will assume that if
$Q_{\beta,l} \subset Q_{\alpha,k}$, then
$\diam(Q_{\alpha,k}) \geq 3^{l-k} \diam(Q_{\beta,l})$ --
to archive this it is enough to replace $Q_{\alpha,k}$
by $Q_{\alpha,mk}$ for $m$ large enough.  

\begin{lem}
Let $Q_{\alpha,k}$ be a family of dyadic cubes on $W_0$.
There exists Borel measurable mapping
$\iota: W_0 \mapsto W$ and constant $C_\iota$ 
such that $\pi\circ \iota = id_{W_0}$
and $\diam(\iota(Q_{\alpha,k})) \leq C_\iota\diam(Q_{\alpha,k}).$
\end{lem}
{\it Proof}: It is enough to construct $\iota$ on
$S = \bigcup_{\alpha}Q_{\alpha,l}$ for some $l$ in such a way that $\iota$
has continuous extension $\iota_\alpha$ to the closure of
each $Q_{\alpha,l}$.
Namely, we can then order linearly the set $I_{l}$ and
extend $\iota$ to whole $W_0$ by taking $\iota(x) = \iota_\alpha(x)$
where $\alpha$ is the smallest one such that $x$ belongs to the
closure of $Q_{\alpha,l}$.  Note that for given $x$ there is only finitely
many such $\alpha$, so $\iota$ is well defined and Borel measurable.
Condition that $\pi(\iota(x)) = x$ is preserved when we take
continuous extension or glue function from pieces.  Also
condition on diameters is preserved: if $x\in Q_{\beta,k}$ for $k <l$
and $x$ is in closure of $Q_{\alpha,l}$ then $Q_{\alpha,l} \subset Q_{\beta,k}$.
So, in the following we will work on $S$ (replace $Q_{\alpha,k}$
by $S\cap Q_{\alpha,k}$) and consequently we will have 
$Q_{\alpha,k} = \bigcup_{Q_{\beta, k+1} \subset Q_{\alpha,k}} Q_{\beta, k+1}$.

Locally we can use smooth $\iota$. So we may assume that
for some $l$ there is mapping $\iota_l$ such that 
$\pi\circ \iota_l = id_{S}$
and for all $k \geq l$ condition on diameters hold.

We choose a maximal family $K \subset I_l$ such that if 
$\alpha_1,\alpha_2 \in K$, $\alpha_1 \ne \alpha_2$, 
$Q_{\alpha_1,l} \subset Q_{\beta_1,k}$,
$Q_{\alpha_2,l} \subset Q_{\beta_2,k}$, then $\beta_1 \ne \beta_2$.
Next, for each $Q_{\alpha,k}$ with $k<l$ we choose $\gamma(\alpha,k)$
such that $Q_{\gamma(\alpha,k), k+1} \subset Q_{\alpha,k}$.  If
$\alpha_0 \in K$ and
$Q_{\alpha_0,l} \subset Q_{\beta,k+1} \subset Q_{\alpha,k}$, then
we choose $\beta$ as $\gamma(\alpha,k)$, otherwise we choose
arbitrarily.  Now, given $\iota_{k+1}$ we want to build $\iota_{k}$.
We do this on each $Q_{\beta, k+1} \subset Q_{\alpha,k}$ separately.
On $Q_{\gamma(\alpha,k), k+1}$ we put $\iota_{k}=\iota_{k+1}$.
If $Q_{\beta, k+1} \subset Q_{\alpha,k}$ and $\beta \ne \gamma(\alpha,k)$
then put $x_0 = x_{Q_{\gamma(\alpha,k), k+1}}$, 
$x_\beta = x_{Q_{\beta,k+1}}$ and 
we choose $y_\beta \in W$ such that 
$d(y_\beta, \iota_{k+1}(x_0)) = d(x_\beta, x_0)$ and
$\pi(y_\beta) = x_\beta$.
Next, put
$\iota_{k}(x) = \iota_{k+1}(x)(\iota_{k+1}(x_\beta))^{-1}y_\beta$.

Observe that sequence $i_k$ is convergent: because $K$ is maximal
for each $Q_{\alpha,k}$ there exists $\beta_0 \in K$ and
$Q_{\beta_1, m}$ such that
$Q_{\beta_0, l} \subset Q_{\beta_1, m}$ and
$Q_{\alpha,k} \subset Q_{\beta_1, m}$.  By our choice of $\gamma(\alpha,k)$
we will have $\iota_{n} = \iota_{m}$ on $Q_{\beta_1, m}$ for all $n< m$.

Finally, we need to check diameter condition.  By our definition
$$\diam(\iota(Q_{\alpha,k})) 
\leq \max_{\beta_1,\beta_2 \in J_{\alpha}}
(d(y_{\beta_1},y_{\beta_2}) + \diam(\iota(Q_{\beta_1, k+1}))
+\diam(\iota(Q_{\beta_2, k+1})))$$
$$
\leq 2\max_{\beta\in J_{\alpha}}
(d(y_{\beta}, y_0)+\diam(\iota(Q_{\beta, k+1})))$$
$$
\leq 2\diam(Q_{\alpha,k}) + 
2\max_{\beta\in J_{\alpha}}\diam(\iota(Q_{\beta, k+1})).$$
Inductively for $k+j \leq l$,
$$
\diam(\iota(Q_{\alpha,k})) \leq 2\sum_{i=0}^{j-1}2^i\max_{\beta:Q_{\beta,k+i} \subset Q_{\alpha,k}} \diam(Q_{\beta,k+i}) 
+ 2^{j}\max_{\beta:Q_{\beta,k+j} \subset Q_{\alpha,k}}
 \diam(\iota(Q_{\beta,k+j})).$$
Choosing $j$ such that $k+j = l$ and using estimate 
$\diam(Q_{\beta,k+j}) \leq 3^{-j} \diam(Q_{\alpha,k})$ we get
$$\diam(\iota(Q_{\beta,k+j})) = \diam(\iota_l(Q_{\beta,l})) \leq
C\diam(Q_{\beta,l}) \leq C3^{k-l}\diam(Q_{\alpha,k})$$
so
$$
\diam(\iota(Q_{\alpha,k})) 
\leq  2\sum_{i=0}^{\infty}(2/3)^i\diam(Q_{\alpha,k})
+C(2/3)^{k-l}\diam(Q_{\alpha,k}) \leq C_1\diam(Q_{\alpha,k}).$$
$\diamond$

{\it Proof of \ref{basecase}}
We are going to construct a doubling family of doubling sets in $G$.
We fix a family of dyadic cubes in $W_0$.

With each dyadic cube $Q$ we associate metric $d_{Q,0}$ on
$N$ by the formula:
$$
d_{Q,0}(n_1, n_2) =
e^{-Mr_Q}d_N(y_Qn_1y_Q^{-1}, y_Qn_2y_Q^{-1})$$
where $r_Q$ is the diameter of $Q$, $y_Q = \iota(x_Q)$,
$x_Q$ is the centerpoint
of $Q$ and $M$ is a large enough constant (to be specified later).

Next, for each $Q$ with $r_Q \geq 1$ we choose a doubling chain of metrics
$d_{Q,j}$, $j=0\dots k_Q$, such that $d_{Q,0}$ is as above
and $d_{Q,k_Q} = d_{S,0}$ where $S$ is the smallest dyadic
cube strictly containing $Q$.  To use \ref{doublechain} we need to
find $m \geq 2$ such that $m^2d_{S,0}\geq d_{Q,0} \geq m d_{S,0}$.
We have:
$$d_{S,0}(n_1, n_2) = e^{-Mr_S}d_N(y_Sn_1y_S^{-1}, y_Sn_2y_S^{-1})
= e^{-Mr_S}d_1(y_Qn_1y_Q^{-1}, y_Qn_2y_Q^{-1})$$
where $d_1(n_1, n_2) = d_N(zn_1z^{-1}, zn_2z^{-1})$ and 
$z = y_Sy_Q^{-1} = \iota(y_S)(\iota(y_Q))^{-1}$.
Since $x_S, x_Q \in S$,
$d_W(z, e) = d_W(y_S, y_Q) \leq \diam(\iota(S)) \leq C_1 r_S$
and by \ref{ad-metric}
$$
d_1(n_1, n_2) \leq \exp(C_2r_S)d_N(n_1, n_2),
$$
$$
d_{S,0}(n_1, n_2) \leq e^{-Mr_S}e^{C_2r_S}d_N(y_Qn_1y_Q^{-1}, y_Qn_2y_Q^{-1})
= e^{-Mr_S}e^{C_2r_S}e^{Mr_Q}d_{Q,0}(n_1, n2)
$$
so putting $m=\exp(Mr_S-Mr_Q-C_2r_S)$ we have
$$
md_{S,0}\leq d_{Q,0}.$$
Similarly, there is $C_3$ such that
$$
d_{Q,0} \leq e^{Mr_S}e^{C_3r_S}e^{-Mr_Q}d_{S,0}.
$$
Note that $1 \leq r_Q \leq r_S/3$.  Also $C_2, C_3 \geq 0$.
If we choose $M \geq (3/2)(2C_2 + C_3+1)$,
then
$$Mr_S -Mr_Q - C_2r_S \geq (2/3)Mr_S - C_2r_S
\geq (2C_2+C_3+1)r_S - C_2r_S$$
$$ = (C_2+C_3+1)r_S\geq 3,$$
so $m \geq \exp(3) \geq 2$.
Also
$$Mr_S-Mr_Q \geq (2/3)Mr_S \geq (2C_2 + C_3)r_S$$
so
$$Mr_S-Mr_Q-2C_2r_S-C_3r_S \geq 0$$
and
$$2(Mr_S-Mr_Q-C_2r_S) = Mr_S-Mr_Q+C_3r_S + (Mr_S-Mr_Q -2C_2r_S-C_3r_S)$$
$$
\geq Mr_S-Mr_Q+C_3r_S.$$
Consequently
$$
d_{Q,0} \leq e^{Mr_S}e^{C_3r_S}e^{-Mr_Q}d_{S,0} \leq m^2d_{S,0}$$
which ends verification of assumptions of \ref{doublechain}.

Now, let ${\cal A}_Q$ be a family of subsets of $N$
defined as follows. If $r_Q <1$, then ${\cal A}_Q$
consists of all balls in $N$ of radius $r_Q$ with
respect to $d_{Q,0}$ metric.  If $r_Q\geq 1$, then
${\cal A}_Q$ consists of all balls in $N$ with
radius $1$ with respect to some $d_{Q,j}$, $j=0,\dots,k_{Q}$.

Finally, we define ${\cal A}$ to be a family of all
sets of form $\iota(Q)R$ where $R\in {\cal A}_Q$.

Note that for $S\in {\cal A}$ diameter $S$ is bounded by
a constant times $r_Q$.  This is obvious if $r_Q \leq 1$,
otherwise we compute:
$$
Sy_Q^{-1} = \iota(Q)Ry_Q^{-1} = \iota(Q)y_Q^{-1}y_QRy_Q^{-1},$$
$$
\iota(Q)y_Q^{-1} \subset B_G(e, C_1r_Q),$$
$$
y_QRy_Q^{-1}\subset B_N(y_Qzy_Q^{-1}, e^{C_2Mr})
\subset B_G(y_Qzy_Q^{-1}, C_3r)$$
where $R = B_{d_{Q,j}}(z,1)$,
so
$$
Sy_Q^{-1} \subset B_G(e, C_1r)B_G(y_Qzy_Q^{-1}, C_3r)
= B_G(y_Qzy_Q^{-1}, C_4r)$$

Also, if $r_Q \leq 1$ then $\iota(Q)R$ is comparable to a ball,
so it is a doubling set.  Otherwise:
$$
\mu(S) = \lambda(Q)\nu(R),$$
$$B_G(e, r_Q)S \subset B_G(e, Cr_Q)y_QR 
\subset \chi(B_W(e, Cr_Q))B_N(e, e^{Cr})y_QR$$
$$
\subset \chi(B_W(e, Cr_Q))y_QB_{d_{Q,j}}(e,1)B_{d_{Q,j}}(z,1)$$
$$
\subset \chi(B_W(e, Cr_Q))y_QB_{d_{Q,j}}(z,2)$$
so
$$
\mu(B_G(e, r_Q)S) \leq \lambda(B_W(e, Cr_Q)y_Q)\nu(B_{d_{Q,j}}(z,2))
= \lambda(B_W(e, Cr_Q))\nu(B_{d_{Q,j}}(z,2))$$
$$
\leq C\lambda(Q)\nu(B_{d_{Q,j}}(z,1)) =
C\lambda(Q)\nu(R)$$
so ${\cal A}$ consists of doubling sets.

It remains to prove that ${\cal A}$ is a doubling family.
Note that second condition of the definition of doubling
family holds, namely, if $S = \iota(Q)R$ and
the metric corresponding to $S$ is not the smallest
one in the chain, than we can take the next metric in
the chain and obtain a ball $T$ such that $R\subset T$,
$\iota(Q)T \in {\cal A}$, 
$2\nu(R)\leq \nu(T)\leq C\nu(R)$.  If the metric corresponding
to $S$ is the smallest in the chain (in particular if $r_Q \leq 1$),
than we can replace the cube $Q$ by bigger one.  In both
cases by \ref{unidouble} we have control of volume.

Again, first condition of the definition of doubling
family for $r_Q\leq 1$ is easy. To get it for $r_Q > 1$, notice
that it is enough to look at
$$
\tilde S_2 =
\bigcup_{S_1 \in {\cal A}, S_1\cap S_2 \ne \emptyset, \mu(S_1)\leq \mu(S_2)}
S_1.$$
Namely, according to the remark after definition we may enlarge $S_2$ to
deduce original condition.

Consider
$S_1,S_2\in {\cal A}$ such that $S_1\cap S_2 \ne \emptyset$.
Then $Q_1\cap Q_2 \ne \emptyset$ and since $Q_1$ and $Q_2$
are dyadic cubes either $Q_1 \subset Q_2$ or $Q_2 \subset Q_1$.
Since $\mu(S_1) \leq \mu(S_2)$, $Q_1 \subset Q_2$.
Also, then the metric corresponding to $S_1$ is greater or
equal to the metric corresponding to $S_2$ -- if $Q_1$ is
strictly smaller than $Q_2$, then this follows form our construction,
if $Q_1 = Q_2$, then bigger metric implies smaller volume
of balls.  Consequently, if $S_2 = \iota(Q_2)B_{d_{Q_2,j}}(z,1)$
$$
S_1 \subset \iota(Q_2)B_{d_{Q_2,j}}(z,2)$$
Like previously, if $d_{Q_2,j}$ is not last in the chain we can
take $\tilde S_2 \subset T = \iota(Q_2)B_{d_{Q_2,j+1}}(z,1)$,
if $d_{Q_2,j}$ is last in the chain we enlarge $Q_2$ first and
then take the next metric in the chain.

%We build a tree of metrics such that:
%$$
%e^{Mc^n}d \leq \rho_{n,\alpha, l} \leq e^{8Mc^n}d$$

$\diamond$

\section{Reduction to the base case}

\begin{defn}\label{CZLS} 
We say that the space~$M$ with metric~$d$ and Borel measure~$\mu$ has
the \emph{Calder\'on--Zygmund property on large scales} if there is
a decomposition like in \ref{CZ}, but $r_i>1/2$ and $g$ may be 
unbounded and only satisfies:
$$
\int_{B(x,1)}|g|d\mu \leq C\lambda\mu(B(x,1)).
$$
\end{defn}

\begin{defn}\label{CZSS}
We say that the space~$M$ with metric~$d$ and Borel measure~$\mu$ has
the \emph{Calder\'on--Zygmund property on small scales} if there is
a decomposition like in \ref{CZ}, but only for $f$ which satisfies:
$$
\int_{B(x,1)}|f|d\mu \leq C\lambda\mu(B(x,1)).
$$
\end{defn}

\begin{lem}\label{czeqczls}
If $G$ is a connected Lie group,
$d$ is a right-invariant optimal control metric on $G$ then $G$ 
(with $d$ and left Haar measure) satisfies Calder\'on--Zygmund property
if and only if it satisfies Calder\'on--Zygmund property on large
scales.
\end{lem}

{\it Proof}: It is well known that optimal control metric satisfies
doubling property for balls of bounded radius: for each $R$ there
is a $C$ such that
$$
\mu(B(x,2r)) \leq C\mu(B(x,r))
$$
for $r \leq R$. This doubling property implies Calder\'on--Zygmund
property on small scales. Together with Calder\'on--Zygmund
property on large scales we get (full) Calder\'on--Zygmund
property.

To get the opposite implication, we fix $\lambda$ and
apply Calder\'on--Zygmund
property to $f$. We will denote by $C_{CZ}$ the constant from
the Calder\'on--Zygmund property. Without loss of generality we
may pretend that $C_{CZ}\geq 2$.

Decomposition given by Calder\'on--Zygmund property may fail
conditions of
Calder\'on--Zygmund property on large scales
because some of $r_i$ may be smaller than or equal to $1/2$, so
we want to correct this. Choose a maximal collection of
non-intersecting balls $\{B(x_\alpha, C_{CZ}/2)\}$ of radius $C_{CZ}/2$.
%(where
%$C_{CZ}$ is the constant from the Calder\'on--Zygmund property). 
If $r_i\leq 1/2$, then $Q_i\subset B(x_i, C_{CZ}/2)$. Since the collection
$\{B(x_\alpha, C_{CZ}/2)\}$ is maximal, there is an $\alpha$ such that
$B(x_i, C_{CZ}/2)$ and $B(x_\alpha, C_{CZ}/2)$ intersect, so 
$$B(x_i, C_{CZ}/2)\subset B(x_\alpha, C_{CZ}).$$
For each $i$ with $r_i \leq 1/2$ choose one $\alpha_i$ as above. 
Put
$$
E_{\beta} = \bigcup_{\alpha_i = \beta} Q_i,
$$
$$
h_{\beta} = \sum_{\alpha_i = \beta} f_i
$$
We have
$$
\int h_{\beta} = 0,
$$
$$
E_{\beta} \subset B(x_\beta, C_{CZ}).
$$
%$h_{\beta} = 0$ outside $E_{\beta}$, so also $h_{\beta} = 0$ outside
%$B(x_\beta, C)$.
Let $A$ be set of all $\beta$ such that
$$
\int |h_{\beta}| \leq C_{CZ}\lambda \mu(B(x_{\beta}, C_{CZ})). 
$$
Put
$$
\tilde g =  g + \sum_{\beta\in A} h_{\beta}.
$$
Now we construct a new family of functions taking those $f_i$ 
for which $r_i > 1/2$ and we add to it all $h_\beta$ with
$\beta\notin A$. We associate set $E_{\beta}$ and radius $1$ with 
each  $h_{\beta}$ for $\beta\notin A$. We claim that this new family
satisfies all conditions required by the Calder\'on--Zygmund property
on large scales (with a new constant).
Namely, we have:
$$
f = \tilde g + \sum_{r_i>1/2}f_i + \sum_{\beta\notin A} h_{\beta},
$$
$$
\sum \|h_{\beta}\|_{L^1} \leq \sum \|f_i\|_{L^1} \leq C_{CZ}\|f\|_{L^1}.
$$
Next
$$
\mu(E_\beta^{*}) \leq \mu(B(x_\beta, C_{CZ}+1))
\leq \mu(B(x_\beta, 2C_{CZ}))
\leq C_1\mu(B(x_\beta, C_{CZ}))
$$
the last inequality due to doubling property for radii up to $C_{CZ}$.
Then
$$
\sum_{\beta\notin A}\mu(E_\beta^{*}) 
\leq C_1\sum_{\beta\notin A}\mu(B(x_\beta, C_{CZ}))
\leq {C_1\over C_{CZ}\lambda}\sum_{\beta\notin A} \|h_{\beta}\|_{L^1}
\leq C_1{\|f\|_{L^1}\over \lambda}
$$
It remains to show that averages of $\tilde g$ over balls of radius
$1$ are bounded. Let $I_x = \{\beta\in A: B(x_\beta, C_{CZ})\cap B(x, 1) 
\ne \emptyset \}$. The cardinality of $I_x$ is bounded by a constant
$C_2$, since the balls $B(x_\beta, C_{CZ}/2)$ are disjoint
and we have
doubling property for radii up to $2C_{CZ}$. Now,
$$
\int_{B(x,1)} |\tilde g| 
\leq  \int_{B(x,1)} |g| + \sum_{\beta \in I_x}\int |h_\beta|
$$
$$
\leq C_{CZ}\lambda\mu(B(x,1)) +
\sum_{\beta \in I_x} C_{CZ}\lambda\mu(B(x_\beta, C_{CZ}))
$$
$$
\leq C_{CZ}\lambda\mu(B(x,1)) + C_2C_{CZ}\lambda\mu(B(x,2C_{CZ}))
$$
$$
\leq C_3\lambda\mu(B(x,1)).
$$
again by doubling for radii up to $2C_{CZ}$.
$\diamond$

\begin{cor}
If a connected Lie group satisfies 
Calder\'on--Zygmund property for one 
optimal control metric, then is satisfies 
Calder\'on--Zygmund property for any other
optimal control metric.
\end{cor}

Namely, by \ref{czeqczls} only behavior at large scales matters,
and at large scales any two optimal control metric are equivalent.
Hence, in the sequel we will usually say that a Lie group satisfies (or not)
Calder\'on--Zygmund property, without mentioning metric.

\begin{lem}\label{compeq}
Let $G_1$ and $G_2$ be connected Lie groups and let $K\subset G_1$ be
a compact Lie group. Assume that $G_2\subset G_1$ and $G_1 = KG_2$ 
or that $K$ is normal and $G_2 = G_1/K$. Then $G_1$ 
satisfies 
Calder\'on--Zygmund property if and only if $G_2$ 
satisfies 
Calder\'on--Zygmund property.
\end{lem}

{\it Proof:} By \ref{czeqczls} we may consider 
Calder\'on--Zygmund property on large scales.
In more detail, consider case when $K$ is normal.
Let $\pi$ be canonical projection from $G_1$ to the quotient.
By choosing appropriate metrics we may assume that metric
balls in $G_2$ are images of balls in $G_1$, that is
$\pi(B(x, r)) = B(\pi(x), r)$.  Consider
$f \in L^1(G_2)$.  Put $h = f \circ \pi$.
Then $\|h\|_{L^1} = \|f\|_{L^1}$.  We use Calder\'on--Zygmund
property on large scales in $G_1$ and obtain $h = \sum h_i + s$.
Let $f_i(x) = \int_{y\in xK}h_i(y)$ and similarly $g(x) = \int_{y\in xK}s(y)$.
We have $f = \sum f_i + g$, $\int f_i = 0$, $\|f_i\|_{L^1} \leq \|h_i\|_{L^1}$.
On $G_1$ we have $h_i = 0$ outside $Q_i$ so putting $R_i = \pi(Q_i)$
we have $f_i = 0$ outside $R_i$ and since $Q_i \subset B(x_i, Cr_i)$
we have $R_i \subset \pi(B(x_i, Cr_i)) = B(\pi(x_i), Cr_i)$.  We
have $R_i^{*} = \pi(Q_i^{*})$, and $\mu(R_i^{*}) = \mu(\pi^{-1}(R_i^{*}))
= \mu(KQ_i^{*})$.  In general $KQ_i^{*}$ has bigger measure than
$Q_i^{*}$, so this would cause trouble.  However,
$Q_i^{*} = \{x : d(x, Q_i) < r_i\}$ and
$KQ_i^{*} \subset \{x : d(x, Q_i) < r_i + \diam(K)\}$.  By
changing $r_i$ so that new $r_i + \diam(K)$ is less or equal
to old $r_i$ we get condition on measure of $R_i^{*}$.  Such
change is possible if old $r_i$ is bigger or equal to $\diam(K)$,
namely we divide $r_i$ by $2$.  Of course, to make sure that
$R_i \subset B(x_i, Cr_i)$ after such change of $r_i$ we need
to multiply $C$ by $2$ to compensate.  When $r_i$ is bigger
than $1/2$ and smaller
than a fixed constant we may estimate measure of $R_i^{*}$
using doubling condition on measure of balls.  Next,
$$\int_{B(x, 1)} |g| \leq \int_{d(y, xK) < 1} |s|(y).
$$
$\{y : d(y, xK) < 1\}$ can be covered by a bounded number of
unit balls and by enlarging $C$ we get condition on $g$.
So $G_2$ satisfies Calder\'on--Zygmund property on large scales.

When proving that $G_1$ has Calder\'on--Zygmund property on large scales
using property of $G_2$ we have problem with defining functions
$f_i$.  To avoid this problem note that it is enough to build
decomposition of nonnegative function: by combining decompositions
of positive and negative parts we get decomposition of arbitrary
function with slightly worse constants.  When $f$ is nonnegative
we put $h(x) = \int_{y\in xK}f(y)$ and apply decomposition to
$h$ obtaining $h = \sum f_i + s$.  We put $t(x) = f(x)/h(\pi(x))$
if $h(\pi(x)) > 0$ and $t(x) = 1$ otherwise.  Now, take
$f_i(x) = t(x)h_i(\pi(x))$ and $g(x) = t(x)s(\pi(x))$.
When $h(\pi(x)) > 0$, then
$$
\sum f_i(x) + g(x) = t(x)(\sum h_i(\pi(x)) + s(\pi(x)) =
t(x)h(\pi(x)) = f(x)$$
When $h(\pi(x))$ in similar way $\sum f_i(x) + g(x) = 0$ and
$f(x) = 0$ so again we get equality (modulo null sets).
Note
that $\int_{y\in xK}t(y) = 1$, so $\int f_i = \int h_i = 0$.
Also, $t(y) \geq 0$, so $|g|(y) = t(y)|s(\pi(y))|$.  Consequently
$$
\int_{B(y, 1)}|g|(y) \leq \int_{B(y, 1)K}t(y)|s(\pi(y))| =
\int_{B(\pi(y), 1)}|s|
$$
and
$$
\|f_i\|_{L^1} = \|h_i\|_{L^1}
$$
We get $Q_i$ in $G_1$ as counterimages of sets from $G_2$
and we take as $x_i$ arbitrary counterimage of corresponding
point in $G_2$.  The only nontrivial condition is
$Q_i \subset B(x_i, Cr_i)$ since $Q_i$ in $G_1$ are bigger
than corresponding sets in $G_2$.  However, like in previous
case we have $Q_i \subset B(x_i, Cr_i + \diam(K))$.  Since
$r_i > 1/2$ by enlarging $C$ we can ensure that new $Cr_i$
is bigger than old $Cr_i + \diam(K)$,
so we obtained Calder\'on--Zygmund decomposition on large
scales on $G_1$.

Case when $G_1 = KG_2$ is similar, but more messy.  To map
functions from $G_1$ to $G_2$ we integrate on cosets of $K$.
To map back we use measurable selector from $K/(K\cap G_2)$
into $K$ like in lemma \ref{messel}.
We transfer sets from $G_2$ to $G_1$ multiplying them by $K$.
We go back using measurable selector.  Now, both ways we can enlarge
sets so we need to compensate as above.

$\diamond$

{\bf Remark}.  In similar way like lemma \ref{compeq} we
could transfer doubling family of doubling sets between
$G_1$ and $G_2$.  More precisely, we need to assume that
sets in the family of diameter of order $1$ are comparable
to balls.  We can compensate changing sizes like above
using remark after definition of doubling family.

\begin{lem}\label{sol-simp}
Let $G$ be a simply connected solvable Lie groups with the
Lie algebra ${\bf g}$ and let $W\subset G$
be a Lie subgroup of $G$ corresponding to a Cartan subalgebra ${\bf w}$
of ${\bf g}$ and $N$ be a Lie subgroup of $G$ corresponding to
${\bf n} = [\bf g, \bf g]$.  Then $W$ and $N$ are simply
connected and $G = WN$.
\end{lem}

{\it Proof:}
Since ${\bf w}$
is a Cartan subalgebra of ${\bf g}$ we have
$$
{\bf g} = {\bf w} + {\bf n}.
$$
Put ${\bf w}_0 = {\bf w} \cap {\bf n}$.
Let $v_1,\dots,v_k \in {\bf w}$ span complementary subspace to
${\bf w}_0$.  Put ${\bf v}_i = \lin\{v_i,\dots,v_k\} + {\bf n}$.
Since $[{\bf g}, {\bf g}] \subset {\bf v}_i$ each ${\bf v}_i$ is
an ideal in ${\bf g}$ and conseqently ${\bf v}_{i+1}$ is an
ideal in ${\bf v}_{i}$.  This means that ${\bf v}_{i}$ is
a semidirect product of one dimensional Lie algebra generated
by $v_i$ and ${\bf v}_{i+1}$.  Let $V_i$ be simply connected
Lie group with the Lie algebra ${\bf v}_{i}$.  Since ${\bf v}_{i}$ is
a semidirect product $V_i$ is also a semidirect product.
Consequently, mapping $(t, g) \mapsto \exp(tv_i)g$ from
${\mathbb R}\times V_{i+1}$ into $V_{i}$ is one to one.
Let ${\tilde N} = V_{k+1}$ be  simply connected
Lie group with the Lie algebra ${\bf n} = {\bf v}_{k+1}$.
Simple induction shows that mapping
$(t_1,\dots, t_k, g) \mapsto \exp(t_1v_1)\dots\exp(t_kv_k)g$
from ${\mathbb R}^k\times {\tilde N}$ into $V_1 = G$ is one to one.
Consequently ${\tilde N} = N$ and $N$ is simply connected.
Similarly, we show that restriction of the mapping above to
${\mathbb R}^k\times {\bf w}_0$ is one to one onto $W$,
so also $W$ is simply connected.  Finally, $\exp(t_iv_i) \in W$
so the expression above means that $G = WN$.
$\diamond$

\begin{lem}\label{sol-center}
Let $G$ and $W$ be as in Lemma \ref{sol-simp}.
Then center of $G$ is contained in $W$.
\end{lem}

{\it Proof:}  We will use notaion form \ref{sol-simp}.
First note that that action of $W$ on $G/W$ by inner
authomorphisms of $G$ is equivalent to action on $N/(W\cap N)$.
Let $z$ be in the center of $G$.
Since $z$ is in the center it is invariant under inner
authomorphizms and it leads to a fixpoint of the
action of $W$ on $G/W$.  So it is enough to show that
the only fixpoint of action on $W$ on $G/W$ by inner
authomorphisms of $G$ is $W$.  This is equivalent to
showing that only fixpoint of action of $W$ on
$N/(W\cap N)$ is $(W\cap N)$.  However, $N$ is a simply
connected nilpotent Lie group, so exponential mapping
from ${\bf n}$ into $N$ is one to one.  So it is enough
to show that the only fixpoint of action of $W$ on
${\bf n}/{\bf w}_0$ is ${\bf w}_0$.  This is equivalent
to showing that the only element of ${\bf n}/{\bf w}_0$
anihilated by ${\bf w}$ is ${\bf w}_0$.  But this holds
since ${\bf w}$ is a Cartan subalgebra of ${\bf g}$.
$\diamond$

Now we can reduce the general case to the base case.

{\bf Step} 1. By the structural theory, amenable Lie group $G$ is
a compact extension of a solvable Lie group $R$. Moreover, by 
Levy-Maltsev theorem there is a compact subgroup $K$ of $G$
such that $G=KR$. So by \ref{compeq} the general case reduces to
solvable one.

{\bf Step} 2. $G$ is now solvable. Let $\alg$ be the complexification of
the Lie algebra of $G$
and let $\ala$ be the Cartan subalgebra of $\alg$. Consider root
space decomposition of $\alg$:
$$
\alg = \bigoplus_{\alpha}\alg_{\alpha}
$$
where $\alg_{\alpha}= \{x \in \alg: 
\forall_{y \in \ala} (\ad (y) - \alpha(y))^n x = 0\}$ and $n$ is 
the dimension of $\alg$. We have $\alg_0 = \ala$ and (by )
$$
[\alg_{\alpha},\alg_{\beta}] \subset \alg_{\alpha+\beta}.
$$ 
Hence, for $y\in\ala$ linear mapping of $\alg$ which
multiplies $x \in \alg_{\alpha}$ by imaginary part of $\alpha(y)$
is a derivation of $\alg$.  Such derivations
commute with complex conjugation on $\alg$, so they generate
automorphisms of universal covering $\tilde G$ of $G$.
Resulting automorphisms act as identity on $\exp{\alg_0}$
which by Lemma \ref{sol-center} contains center of $\tilde G$.
Conseqently we can pass to quotient and obtain automorphisms
of $G$.

Obviously, the closure $T$ of all such 
automorphisms is a group isomorphic to a torus (and compact). It is
easy to see that in the semidirect product $T\cross G$ we can
find a subgroup $G_1$ such that all roots of $G_1$ are real.
So, applying \ref{compeq} twice, first to $G$ and $T\cross G$,
then to $T\cross G_1$ and $G_1$ we reduce the solvable case to
the case with all roots real. 

{\bf Step} 3. Now $G$ has all roots real. What matters for us, is
that the Lie algebra $\alg$ of $G$ is exponential: the exponential
mapping from $\alg$ to corresponding simply connected Lie group
$\tilde G$ (which we will identify with universal covering of $G$)
is a diffeomorphism. It follows that the center $Z$ of $\tilde G$ is
an image of the center of $\alg$. Also, the kernel $N$ of the covering
map from $\tilde G$ to $G$ is a central subgroup. So, we may
identify $N$ with a lattice $L$ in the center of $\alg$. Let $K$ be
a subgroup of $G$ corresponding to linear span $V$ of $L$. $K$
is isomorphic to $V/L$, so it is a normal torus. Again using 
\ref{compeq} we may divide $G$ by $K$ reducing the problem to
simply connected groups.

{\bf Step} 4. Now $G$ is an exponential solvable Lie group. Like
in step 2 we fix a Cartan subalgebra and root space decomposition
of the Lie algebra $\alg$ of $G$:
$$
\alg = \bigoplus_{\alpha}\alg_{\alpha}
$$
(since after step 2 all roots are real we can skip complexification).
Let $\aln$ be the subalgebra of $\alg$ generated by all $\alg_{\alpha}$ with
$\alpha \ne 0$.
Let $N$ be the subgroup of $G$ corresponding to $\aln$.
$N\subset [G,G]$ is a normal divisor. Since $G$ is
exponential $N$ is closed.  Let $G_0$ be
the subgroup of $G$ corresponding to $\alg_{0}$. $G_0$ is nilpotent,
hence of polynomial growth. We need to get estimate on the 
distance in $N$.
 Let $d$ be optimal control metric corresponding
to basis $X_1,\dots, X_m, Y_1,\dots, Y_l$ of $\alg$ consisting of
basis $X_1,\dots, X_m$ of $\aln$ extended by vectors $Y_1,\dots, Y_l$
from $\alg_0$ ($d$ is in fact a Riemannian distance). Consider
a curve $\gamma:[0,1]\mapsto G$ joining $e$ with $x\in N$. We 
have
$$\gamma'(s) = \sum_{i=1}^m a_i(s)X_i(\gamma(s)) + 
\sum_{i=1}^{l}a_{m+i}(s)Y_i((\gamma(s))
$$
Let $\gamma_1$ be the solution of the differential equation
$$\gamma_1'(s) =  
\sum_{i=1}^{l}a_{m+i}(s)Y_i((\gamma_1(s))
$$
with initial condition $\gamma_1(0) = e$.
Consider images $\pi\circ\gamma$ and $\pi\circ\gamma_1$ of
$\gamma$ and $\gamma_1$ under quotient projection $\pi:G\mapsto G/N$.
Both $\pi\circ\gamma$ and $\pi\circ\gamma_1$ satisfy the same
differential equation with the same initial condition, so
$$\pi\circ\gamma(s) = \pi\circ\gamma_1(s).$$
Put $\gamma_2(s) = \gamma_1^{-1}(s)\gamma(s)$. By the above
$\pi(\gamma_2(s)) = e $ so $\gamma_2(s)\in N$.
Next
$$
\gamma_2'(s) = (\gamma_1^{-1}\gamma)'(s) = 
dL_{\gamma_1^{-1}(s)}\gamma'(s) + dR_{\gamma(s)}(\gamma_1^{-1})'
$$
$$
= \Ad(\gamma_1^{-1}(s))dR_{\gamma_1^{-1}(s)}\gamma'(s) 
+ dR_{\gamma(s)}(\gamma_1^{-1})'
$$
On the other hand, since $\gamma_2:[0,1]\mapsto N$ there are
(unique) $b_1,\dots, b_m$ such that
$$
\gamma_2'(s) =  \sum_{i=1}^m b_i(s)X_i(\gamma_2(s))
$$
Our vector fields are right invariant, so for each $s$
$$
\sum_{i=1}^m b_i(s)X_i = 
\sum_{i=1}^m a_i(s) \Ad(\gamma_1^{-1}(s)) X_i
+ \sum_{i=1}^{l}a_{m+i}(s) (\Ad(\gamma_1^{-1}(s))-1)Y_i
$$
and
$$
\sum b_i^2(s) \leq e^{Ct} \sum a_i^2(s)$$
where $t$ is length of $\gamma$. Hence, for $x\in N$, 
$d_G(x, e)< t$ implies $d_N(x, e)< te^{Ct}$,
% so for $t\geq 1$, $B_G(e, t)\cap N \subset B_N(e, e^{C_2t})$,
which in turn implies
$1 + d_N(x) \leq \exp((C + 1)d_G(x))$.

{\bf Step} 5. To obtain remaining estimate on distance we first prove that
$$
B_N(e, e^{t}) \subset B_G(e, C(t+1))\cap N
$$
for $C>0$ large enough and $t>0$. 
As in step 4 we may assume here
that
$X_1,\dots, X_m$ is a basis of
$\bigoplus_{\alpha\ne 0}\alg_{\alpha}$,
such that each $X_i$ belongs to some $\alg_{\alpha}$ and that
$Y_1,\dots, Y_l$ is a basis of $\alg_0$ (and $d$ corresponds to
fields $X_1,\dots, X_m$, $Y_1,\dots, Y_l$).

We need a lemma:

\begin{lem}
%If $X_i$, $i=1,\dots, m$ generate the Lie algebra of $N$
Let $X_i$, $i=1,\dots, m$, $N$ and $d_N$ be as above. There
exists a constant $k$ such that
$$
B_N(e, r) \subset \left(\bigcup_{i=1}^{m}\{\exp(sX_i): |s| <r\}\right)^k
$$
\end{lem}
{\it Proof}: Note that we can replace
$N$ by a free nilpotent group of the same step as our original $N$
(the claim is preserved when passing to a quotient).  Hence, we
may assume that $N$ has a one-parameter family of authomorphic
dilations $\delta_s$, $s>0$, such that
$$\delta_s(X_i) = sX_i$$
so
$$
d(e, \delta_s(x)) = sd(e, x)$$
Since $X_i$ generate $N$ we can find $k$ such that the
inclusion holds for $r=1$. Applying $\delta_s$ to both sides
we see that the inclusion holds for all $r>0$.
$\diamond$

Fix $i$. There is $\alpha$ such that $X_i\in \alg_{\alpha}$.
Fix $Z\in \alg_{0}$ such that $\alpha(Z) = 1$
and $b_1, \dots, b_l$ such that $Z=\sum_{i=1}^{l}b_iY_i$.
Let $\gamma(s) = \exp(-3tsZ)$ for $0\leq s \leq 1/3$,
$\gamma(s) = \exp(3(s-1/3)X_i\exp(-tZ)$ for $1/3\leq s \leq 2/3$,
$\gamma(s) = \exp(3t(s-2/3)Z)\exp(X_i)\exp(-tZ)$ for $2/3\leq s \leq 1$.
We have
$$\gamma(1) = \exp(tZ)\exp(X_i)\exp(-tZ) = \exp(\exp(t\ad(Z)X_i) =
\exp(\exp(t)X_i)$$
$$
\gamma'(s) = \sum_{i=1}^{m}a_iX_i(\gamma(s)) 
+ \sum_{i=1}^{l}a_{i+m}Y_i(\gamma(s))$$
where $|a_i|$ are bounded by $C(t+1)$ with $C$ large enough, so 
$$
\exp(\exp(r)X_i) \in B_G(e, C(r+1)).
$$
By the lemma
$$
B_N(e, \exp(t)) \subset B_G(e, C_2(t+1))$$
with (another) $C_2>0$ large enough.
Now, let $1 < d_N(e, x) = e^t$.  By the inclusion above
$d_G(e, x) \leq C_2(t+1)$, so $C_2^{-1}d_G(e, x) - 1 \leq t$,
so
$$
\exp(C_2^{-1}d_G(e, x) - 1) \leq e^t = d_N(e, x) \leq 1 + d_N(e, x).
$$
When $d_G(e, x) \geq 2C_2^{-1}$ we have
$C_2^{-1}d_G(e, x) - 1 \geq (2C_2)^{-1}d_G(e, x)$, so
$$
\exp((2C_2)^{-1}d_G(e, x)) \leq 1 + d_N(e, x).
$$
$d_G(e, x)$ is locally Lispchitz, so there is constant $C_3$ such
that when $x\in N$ and $d_G(e, x) \leq 2C_2^{-1}$ we have
$d_G(e, x) \leq C_3d_N(e, x)$ and
$$
\exp(C_4d_G(e, x)) \leq 1 + d_N(e, x)
$$
with $C_4 = (C_3\exp(2C_2^{-1}))^{-1}$
which finishes our reduction of \ref{amencz} to \ref{basecase}.

\section{Necessary condition for Calder\'on--Zygmund property}

We say that a set $A$ is $r$-doubling iff $|B(r)A| \leq 2|A|$.

\begin{theo} Assume a group $M$ with a right-invariant metric $d$ and
left-invariant measure satisfies Calder\'on--Zygmund property.
If $M$ has infinite volume
then $M$ contains $r$-doubling sets for arbitrarily large $r$.
%If $M$ contains
%arbitrarily small subsets of positive measure, then it contains
%arbitrarily small doubling subsets.
\end{theo}

{\it Proof}: It is enough to find sets $A(r_i)$ such that 
$|B(r_i)A(r_i)|\leq 2|A(r_i)|$ for a sequence of $r_i\rightarrow\infty$.
We will suppose that the inverse inequality
$|B(r)A| > 2|A|$ holds for all large $r$ and will derive
contradiction with Calder\'on--Zygmund property.

If $M$ has infinite volume, we can take
$\lambda$ arbitrarily small. Let $f(x)=1/|B(1)|$ for $d(x,e)<1$
and $0$ otherwise. Of course, $\norm{f}_{L^1}=1$.
Fix a Calder\'on--Zygmund decomposition of $f$ with small
$\lambda$ (to be determined later).
Put 
$$E_k = \bigcup_{2^k<r_i\leq 2^{k+1}} Q_i,$$
$$F_k = B(1)\cup \bigcup_{j\geq k} E_j,$$
$$h_k = \sum_{2^k<r_i} f_i.$$
Let 
$\phi_k(x)=1$ for $d(x,F_k)>2^{k-1}$ and 
$\phi_k(x)=2^{-k+1}d(x,F_k)$ otherwise.
 
Let $C$ be the constant from the definition of Calder\'on--Zygmund property.
Fix $l>4$ such that $2^{l-3}>C^2$.  Assume that
$|B(2^k+1)|\leq C/\lambda$ and
$2^{k-1}>rl$.
If $r$ is large enough, then there are no $r$ doubling sets
and
$$|B(r)^{j+1}B(2^{k-1})F_k| \geq 2|B(r)^{j}B(2^{k-1})F_k|.$$
By induction
$$|B(2^{k-1})F_k|\leq 2^{-l}|B(r)^{l}B(2^{k-1})F_k|.$$
Next
$$|B(2^{k-1})F_k|\leq 2^{-l}|B(r)^{l}B(2^{k-1})F_k|\leq 2^{-l}|B(2^k)F_k|$$
$$\leq 2^{-l}(|B(2^k+1)|+\sum \mu(Q^{*}_i))
\leq 2^{-l+1}C\norm{f}_{L^1}/\lambda
\leq 1/(4C\lambda)$$
so
$$\int_{B(2^{k-1})F_k} |g| \leq C\lambda|B(2^{k-1})F_k| \leq 1/4.$$
Also
$$\int g = \int f -\sum \int f_i = \int f = 1$$
so
$$\int \phi_kg \geq \int g - \int_{B(2^{k-1})F_k} |g| \geq 1 - 1/4 = 3/4.$$
If $r_i\leq 2^{k-l}$ then
$$\int \phi_kf_i\leq 2^{-k+1}r_i\norm{f_i}_{L^1}\leq 2^{-l+1}\norm{f_1}_{L^1}$$
so
$$
\sum_{r_i\leq 2^{k-l}}\int \phi_kf_i\leq 2^{-l+1}\sum\norm{f_i}_{L^1}
\leq 2^{-l+1} C\norm{f_i}_{L^1} \leq \norm{f}_{L^1}/4 = 1/4.$$
Next
$$\phi_kf = 0,$$
$$\phi_kh_k = 0,$$
so
$$
-\sum_{2^{k-l}<r_i\leq 2^{k}}\int \phi_kf_i = 
\int \phi_kg - \int \phi_kf - \int \phi_kh_k 
- \sum_{r_i\leq 2^{k-l}}\int \phi_kf_i$$
$$\geq 3/4 - 1/4 = 1/2 $$
Now, if $m>2C$, and $\lambda$ is small enough so the inequalities above
hold for $k_j = k_0 +jl$, $j=0, \dots, (m-1)$ then
$$\sum_{2^{k_0+(j-1)l}<r_i\leq 2^{k_0+jl}}\norm{f_i}_{L^1}\geq 
-\sum_{2^{k_0+(j-1)l}<r_i\leq 2^{k_0+jl}}\int \phi_{k_j}f_i \geq 1/2$$
so
$$\sum \norm{f_i}_{L^1} \geq \sum_{j=0}^{m-1}
\sum_{2^{k_0+(j-1)l}<r_i\leq 2^{k_0+jl}}\norm{f_i}_{L^1}
\geq \sum_{j=0}^{m-1}1/2 = m/2 > C\norm{f}_{L^1}$$
which gives a contradiction.

%Reasoning to prove existence of arbitrarily small doubling sets is similar.
$\diamond$

\begin{cor}
If a locally compact group G with a right-invariant geometric metric $d$ and
the left Haar measure satisfies Calder\'on--Zygmund property, then
$G$ is amenable.
\end{cor}

%{\it Proof}: Note that non-amenable group has infinite volume.

\end{document}